\documentclass[12pt,leqno]{amsart}
\usepackage[dvips]{graphics}
\usepackage{amssymb}

\numberwithin{equation}{section}

\newtheorem{thm}{Theorem}[section]

\newtheorem*{thm*}{Theorem}
\newtheorem*{thmmain*}{MAIN THEOREM}
\newtheorem{lem}{Lemma}[section]
\newtheorem{cor}[lem]{Corollary}
\newtheorem{prop}[lem]{Proposition}

\newtheorem{defn}[lem]{Definition}
\theoremstyle{remark}
\newtheorem{rem}{Remark}[section]
\newtheorem{ex}[rem]{Example}

\newcommand{\cref}[1]{Corollary~\ref{#1}}

\newcommand{\R}{\mathbb{R}}

\begin{document}
\title{Affine functions on $CAT(\kappa)$-spaces}
\author{Alexander Lytchak \ \&
Viktor Schroeder}


\maketitle
\renewcommand{\theequation}{\arabic{section}.\arabic{equation}}
\pagenumbering{arabic}

\section{introduction}
This paper is devoted to the structure of singular 
metric spaces 
admitting  
affine functions. Since we are dealing with quite general 
spaces it 
is reasonable
to require the maps to be Lipschitz.

\begin{defn}
A Lipschitz map $f:X \to Y$ between geodesic 
metric spaces is 
called affine,
if it maps each geodesic to a geodesic 
parametrized proportional
to the arclength.
In the case $Y=\R$ we call $f$ an affine function.
\end{defn}

The easiest example of an affine map occurs
in the situation that $X$ splits
as $X' \times Y$ and $f$ is the projection 
$p:X'\times Y\to Y$.
In the case $Y=\R$ we get affine functions. In
 \cite{AB}  situations are studied, in which 
the existence of an affine function 
$f: X \to \R$ already implies the existence of a Euclidean
de Rham factor. To obtain these results
one has to assume that the space is geodesically
complete. Without this assumption one cannot expect the
existence of a splitting. The best one can hope for is 
the existence
of an isometric embedding of
$X$ into a product with a flat factor.
Indeed our main result is

\begin{thm} \label{mainthm}
Let $X$ be a $CAT(\kappa )$ space. Then there is a canonically 
defined isometric 
embedding $i:X \to Y \times H$, where $Y$ is a 
geodesic metric
space and
$H$ is a Hilbert space. 
Every affine function
$f:X \to \R$ factors as
$f = \hat{f} \circ p_H \circ i$ where
$p_H$ is the projection onto $H$ and $\hat{f}:H\to\R$ 
is an affine function.
Moreover each isometry of $X$ determines an isometry 
of $Y$ and of $H$.
Finally the completion of $Y$ is $CAT(0)$ if $X$ is $CAT(0)$.
\end{thm}

\begin{rem}
In the case that $X$ is a Hadamard space 
and the affine function is a Busemann
function this result was shown in \cite{AdB}. 
This was the motivation for
our work.
\end{rem}

\begin{rem}
If we assume (instead of the
$CAT(\kappa)$ condition) that
$X$ is an Alexandrov space with lower curvature bound
and finite dimension (maybe with nonempty boundary)
a corresponding theorem can be proved in 
essentially the same way.
This generalizes results of Alexander and Bishop \cite{AB}.
\end{rem}

Without a curvature assumption a corresponding result is
no longer true:

\begin{ex}
 Let $X$ and $Y$ be geodesic spaces, 
$||\cdot ||$ a strongly
convex norm on a two dimensional vector space. Let 
$Z= X\times _{||\cdot ||} Y$ be the 
non-standard metric product in the sense
of \cite{BFS}. Then
the projections $p:Z\to X$ and $Z\to Y$ are affine. 
In particular if $Y$ is 
a strongly convex Banach space one gets many affine 
functions on $Z$.
Moreover convex subsets of non-standard products admit 
affine functions. Such a space does not necessarily admit
a non-trivial isometric embedding into a space with a direct
Euclidean factor.
\end{ex}

The next example 
describes a more complicated space with a non-trivial 
affine function which doesnot admit an embedding even into a 
nonstandard product.

\begin{ex}
Let $B_1$ and $B_2$ be two Banach spaces with strongly 
convex and smooth 
norms. Let $v_i$ be a unit vector in $B_i$. Denote by 
$\gamma _i$ the
line $\gamma _i (t) = tv_i$ and by 
$f_i:B_i \to \mathbb {R}$ the Busemann
function  of $\gamma _i$. By identifying 
$\gamma _1$ and $\gamma _2$
we glue $B_1$ and $B_2$ to a space $X$. 
Observe now that the function 
$f:X\to \mathbb {R}$ that arises from 
$f_1$ and $f_2$ is affine.
\end{ex}

For general metric spaces it is not clear which implications 
the existence of an affine function has.
Under the additional assumption that the affine
functions separate the points in $X$ one can prove that
$X$ is isometric to a convex subset of a Banach space.

 All the proofs in \cite{IN},\cite{Ma1},\cite{Ma2} and \cite{AB} 
have in common
that the non-Euclidean factor can be recognized as a 
convex subset of 
$X$.  
Our proof is quite different
and the outline of the argument is
as follows:
Let
$\mathcal A$
be the space of affine functions on
$X$
modulo the constant functions.
If
$X$
is 
$CAT(\kappa)$,
then
$\mathcal A$
and its dual space
$H:=\mathcal A^*$
are Hilbert spaces (section 4).
There is a naturally defined
{\em evaluation map}
$F: X \to H$.
In section 5 we prove that the function
$\tilde d:X\times X \to [0,\infty)$,
$\tilde d(y,z) = \sqrt {d(y,z) ^2 - ||F(y) -F(z)|| ^2}$
defines a pseudometric on 
$X$.
Let
$Y = X/{\tilde d}$
be the corresponding metric space.
We finally show that
$i: X \to Y\times H$,
$x \mapsto ([x],F(x))$
satisfies the properties of
Theorem \ref{mainthm}.

\begin{rem}
We note that in general the factor
$Y$
cannot be embedded
isometrically into 
$X$.
This makes it difficult to
obtain geometric properties
of
$Y$.
We do not know, if 
the 
$CAT(\kappa)$ property
of $X$ implies 
$CAT(\kappa)$
for $Y$.
In the special case
$\kappa =0$ we can however
prove this.
\end{rem}


\section{Preliminaries}

By $d$ we will denote the distance in metric spaces 
without an extra
reference to the space.  A {\it pseudo metric} is a 
metric for which
the distance between different points may be zero. It 
defines a unique metric 
space. 

 A {\it geodesic} in a metric space is a
length minimizing curve parametrized
propotionally to arclength. A metric space is 
{\it geodesic} if all pair of points are 
connected by a geodesic. A subspace of a geodesic space is 
{\it convex} if 
it is geodesic with respect to the induced metric. 
A $CAT(\kappa )$ space
is a complete  geodesic metric space in which 
triangles are not thicker
than in the space of constant curvature $\kappa$. 
We refer to  \cite{BH} for more detailed discussion
of these spaces.

A map $f:X \to Y$ is called {\it $L$-Lipschitz} 
if $d(f(x),f(z)) \leq L d(x,z)$.
The smallest $L$ is called the optimal Lipschitz constant.
For a Lipschitz function $f: X \to R$  we denote 
by $|\nabla _x f|$  the {\it absolut gradient} at $x$ 
which is given by
$\sup \{ 0 , \limsup _{z\to x} \frac {f(z)-f(x)} {d(x,z)} \}$. 
If the space $X$ is geodesic,
the optimal Lipschitz constant is the supremum 
of all absolut gradients.   

Remark that a Lipschitz function 
$f:X\to \R$ is affine iff it is convex
 and concave, i.e. if its restriction to each 
geodesic is convex and concave.
For a convex (in particular for an affine) function 
$f$ the absolut 
gradient $|\nabla _x f| $ is semi-continuous in 
$x$ (compare \cite{P}).


\section{ Affine functions on general spaces}

 Let $X$ be an arbitrary geodesic metric space. 
The set of all
affine functions on $X$ is a vector space and 
will be denoted by
$\tilde {\mathcal A }(X)$.  It always contains
the one-dimensional subspace $Const (X)$ of constant
functions.  For each point $x \in X$ the space 
$\mathcal A _x$ of
all affine functions vanishing at $x$ is a complement
of $Const$ in $\tilde {\mathcal A} (X)$. By $\mathcal A (X)$ 
or simply $\mathcal A$
we will denote
the quotient vector space 
$\tilde {\mathcal A} (X) / Const (X)$. 
For an affine  function
$f:X \to \R$ we denote with
$[f] \in \mathcal A$ the corresponding element of
$\mathcal A$.
 The best Lipschitz constant defines a norm on the 
space $\mathcal A$. 
Equipped with this norm $\mathcal A $ 
is a normed vector space.
It is complete (even if $X$ is not complete), hence it is 
a Banach space.

Consider the evaluation map 
$E:X\times X \to \mathcal A ^* $ 
from the product $X \times X$
to the dual space of $\mathcal A$ given by
$E(x,y) ([f]) =f(y) -f(x)$.
We have 
$$E(x,z) ([f]) - E(\bar x, \bar z) ([f])
\leq ||f|| \  (\ d(x,\bar x) +d(z,\bar z)\ ).$$
Moreover the map $E$ is strongly 
affine in the sense that it maps geodesics to affine lines of
the Banach space $\mathcal A^*$. 
Observe that $E(x,z)=0$ iff the points $x$ and $z$ cannot
be separated by an affine map on $X$. 

By $E_x :X \to \mathcal A ^*$ we 
denote the restriction $E_x (z) ([f]) =f(z)-f(x)$.
We have $E_y =E_x + E(x,y)$.


\section{Affine functions on $CAT(\kappa)$ spaces}

Let $X$ be a 
$CAT(\kappa)$ space and
$f:X \to \R$ affine.
For
$x\in X$ let
$C_x = CS_x$ be the tangent cone at the point 
$x \in X$ which is the cone
over the space of directions $S_x$.
Then $f$ induces a homogeneous
affine function ( the {\it directional derivative})
$D_x f: C_x \to \R$ (compare \cite{K}).
The absolute gradient
$\mid \nabla_x f \mid$
is equal to
$\sup \{D_xf(v) \mid v\in C_x, d(0,v)=1\}$.
The function
$D_xf$ inherits the Lipschitz constant from $f$.

The following splitting result is basic:

\begin{lem} \label{splitting}
Let $X$ be a $CAT(0)$ space. 
Let $f:X\to \R$ be an 
affine function. Assume that for some line $\gamma$ 
in $X$ we have
$(f\circ \gamma ) ' = ||f|| >0$. 
Then $X$ splits as  $X=Z\times \R$ and
$f$ is given by $f(z,t) =||f||t$.
\end{lem}

\begin{proof}
We may assume that $||f||=1$ and $f(\gamma (0)) =0$.
Let $x\in X$ be arbitrary. For the
rays $\gamma _x ^+$ and $\gamma _x ^-$ starting at $x$ and 
asymptotic to 
$\gamma ^+$ resp. $\gamma ^-$ we immediatly obtain
$(f\circ \gamma_x ^+)' =1$ and 
$(f\circ \gamma_x ^-)' = -1$. 
Therefore
$|f(\gamma_x^+(1)) - f(\gamma_x^-(1))| = 2$.
Since $f$ is $1$-Lipschitz we deduce
that
$d(\gamma_x^+(1),\gamma_x^-(1)) = 2$
and hence
the concatenation of $\gamma_x^+$ and $\gamma_x^-$ is
a line $\gamma_x$ which is parallel
to $\gamma$. 
Therefore  through each point $x\in X$ there is a line
paralell to $\gamma$ and we may 
apply the well known splitting theorem
(\cite{BH}). Now the last statement is clear too.
 \end{proof}

\begin{prop}
Let $X$ be a $CAT(\kappa)$ space and
$f:X \to \R$  an affine function. 
Assume that $y$ is an inner point of a 
geodesic starting at $x$. 
Then $|\nabla _y f| \geq |\nabla _x f|$.
\end{prop}

\begin{proof}
Let $X$ be a $CAT(\kappa)$ space. We may assume 
that $d(x,y)< \frac{\pi}{3\sqrt{\kappa}} $ and
that for some point $z$ we have
$d(z,y)=d(x,y)= \frac 1 2 d(x,z)$. Moreover
we may assume $f(x)=0$. Let $f(y)=r$.
Let $\eta $ be a geodesic starting at $x$ with 
$a= (f\circ \eta ) ' >0$.
Consider the midpoint $m_t$ of the 
geodesic between $z$ and $\eta (t)$
for small $t$. We have $f(m_t)= \frac {2r + at } 2$. 
On the other hand the $CAT(\kappa)$ assumption implies
$d(y,m_t) \leq \frac t 2 + At^2$ for some fixed 
$A \geq 0$ depending only on $\kappa$. This implies
$$|\nabla _y f| \geq  \frac{f(m_t)-f(y)}{d(m_t,y)}
\geq \frac{\frac{at}{2}}{\frac{t}{2}+At^2} = \frac{a}{1+2At}$$
For $t \to 0$ we obtain 
$|\nabla _y f| \geq a$.
Since $\eta$ is arbitrary we have
$|\nabla _y f| \geq  |\nabla _x f|$.
\end{proof}

We see that  for each affine  function $f$ on $X$, the set 
$X_{\epsilon}$ of all points $x\in X$ such 
that $|\nabla _x f| > ||f|| - \epsilon$ is
open, dense and convex in $X$. 
From the theorem of Baer we obtain:

\begin{cor} \label{baer}
 Let $f_k$ be a sequence of affine function.  Then the set
$X ^0$ of points $x$ such that   
$|\nabla _x (-f_j)|= |\nabla _x f_j| =
|| f _ j||$ for all $j$ is convex and dense in $X$.
\end{cor}

Now we can deduce

\begin{cor} \label{hilbert}
 Let $X$ be a $CAT(\kappa)$ space. Then
 $\mathcal A $ is a Hilbert space.
\end{cor}

\begin{proof}
Let $f,g$ be two affine functions. We have to prove 
$||f+g||^2 + ||f-g||^2 = 2 (||f||^2 + ||g||^2)$.
By Corollary  \ref{baer} there exists
$x\in X$ sucht that
$|\nabla _x h| = ||h||$ for $h=\pm f,\pm g,f+g,f-g$.
For simplicity let
$h':=D_xh: C_x \to \R$.
Then the functions
$h'$ are homogeneous, affine and
$||h'|| =||h||$.
Let $0\in C_x$ be the cone point of
$C_x = CS_x$,
where $S_x$ is the space of directions in $x$.
Since $||\pm f'|| = ||f||$ we have
$|\nabla_0 f'| = ||f||$ and
$|\nabla_0 (-f')| = ||f||$.
Thus there are
$v^+,v^- \in S_x$ such that
$f'(v^+) = -f'(v^-)=||f'||$,
hence
$|f'(v^+)-f'(v^-)| = 2 ||f'||$.
This implies
$d(v^+,v^-) = 2$, where this distance is
measured in the cone
$CS_x$. Thus the concatenation of the two rays
$\gamma^+(s) = s v^+$ and
$\gamma^-(s) = s v^-$ for $s \in [0,\infty)$
is a line in the cone $C_x$, and
$(f'\circ\gamma)' = ||f'||$ along this line.
By Lemma \ref{splitting} the $CAT(0)$ space
$C_x$ splits as
$Z\times \R$ and
$f'(z,t) =||f'|| \cdot t$.
In the same way
$C_x$ can be decomposed
as
$Z' \times \R$ such that
$g(z',s) = ||g'|| \cdot s$.
By the properties of the Euclidean de Rham
factor of a $CAT(0)$ space
(compare \cite{BH} p. 235),
$C_x$ splits as
$Z''\times E$, where $E$ is a one or
twodimensional Euclidean space
and
$f' =\hat{f'}\circ p_E$ ,
$g' =\hat{g'}\circ p_E$ ,
where $p_E$ is the projection onto $E$
and
$\hat{f'},\hat{g'}$ are affine functions
on the Euclidean space.
Thus the equation
$||f'+g'||^2 + ||f'-g'||^2 = 2 (||f'||^2 + ||g'||^2)$
and hence the corresponding equation for
$f,g$ holds.

\end{proof}

We come back to the affine maps
$E_x:X \to \mathcal A^*$
defined in section 3.
In the case that
$\mathcal A^*$ is a Hilbert space,
these maps are normalized in the following sense.

\begin{defn}

Let  $X$ be a geodesic metric space,
$H$ be a Hilbert space and  
$F:X\to H$
an affine map.  We call $F$ normalized, 
if $F$ is $1$-Lipschitz
and for each unit  vector $v\in H$ the affine 
function $F^v:X\to \mathbb{R}$
given by $F^v (x) =\langle F(x),v \rangle$ 
satisfies $||F^v||=1$. 
\end{defn}

\begin{ex}
 Let $H_0 \subset H$ be a  Hilbert subspace. 
Then the orthogonal projection
$p:H\to H_0$ is normalized . 
If $F:X\to H$ is normalized, then so is
the composition $p\circ F$.
\end{ex}

Observe that if $F:X\to H$ is normalized, 
then the linear hull of 
the convex set $C=F(X)$ is dense in $H$. 
By the very definition the canonical
evaluation maps $E_x :X\to \mathcal A^*$ are normalized.

\begin{defn}
 Let $X$ be a $CAT(\kappa)$ space, $H$ a Hilbert space and
$F:X\to H$ a normalized affine map.
We call a point $x\in X$ regular if $C_x$ has the splitting 
$C_x =C_x' \times H_x$, with a Hilbert space $H_x$, such 
that
$D_x F $ is the projection onto $H_x$.
\end{defn}

\begin{cor} \label{regular}
If $H$ is separable,
then the set of regular 
points is convex and dense
in $X$.
\end{cor}

\begin{proof}
Let $e_i$, $i \in \mathbb N$
be a dense subset of the
unit vectors in $H$
and let
$F_i=F^{e_i}$
be the corresponding affine functions.
By Corollary \ref{baer}
the set
$W \subset X$ of points $x$
such that
$|\nabla_x(F_i)| =|\nabla_x(-F_i)| = ||F_i||$
for all $i\in \mathbb N$
is convex and dense.
For $x \in W$ let
$C_x$ be the tangent cone which is a 
$CAT(0)$ space and splits an Euclidean
de Rham factor
$C_x = C_x''\times H'_x$.
By the proof of Corollary \ref{hilbert}
the homogeneous affine function
$D_xF_i:C_x \to \R$
has the form
$D_xF_i(z'',h') = \langle v_i,h'\rangle$,
where
$v_i \in H_x'$
is a unit vector.
Let
$H_x \subset H'_x$ be the 
closure of the span of
the $v_i$
and 
$C_x =C'_x\times H_x$ be the corresponding
splitting where
$C'_x =C''_x\times H_x^{\perp}$.
By construction 
$D_xF$ is the projection onto
$H_x$.
\end{proof}

\section{Proof of Theorem \ref{mainthm}}

The proof of the main theorem is based on the following
fact

\begin{thm}
Let $X$ be a $CAT(\kappa)$ space, $H$ be a Hilbert space and 
$F:X\to H$ be a normalized affine map. 
Then $\tilde d :X\times X \to [0,\infty)$ given by
$\tilde d (y,z) = \sqrt {d(y,z) ^2 - ||F(y) -F(z)|| ^2}$
defines a pseudo metric on $X$.
\end{thm}

\begin{proof}
 
By definition
$\tilde d$ is symmetric and since
$F$ is 1-Lipschitz, $\tilde d$ is
nonnegative.

Since $F$ is affine, we have for each
point $m$ on a geodesic $yz$ that

\begin{equation} \label{inpoint}
\tilde d (y,z)= \tilde d(y,m) + \tilde d (m,z).
\end{equation}

We will prove that $\tilde d$ satisfies the triangle
inequality and therefore defines a pseudometric.
We will first show that the triangle inequality is 
satisfied in the neighborhood of every point.
Consider therefore three points $x,y,z \in X$ with
pairwise distance $< \frac{\pi}{2 \sqrt{\kappa}}$ and
assume that
$\tilde d(x,z) > \tilde d (x,y) + \tilde d (y,z)$.
We may assume that
$F(x) = 0$.
Denote by $H_0$ the linear hull
of $F(y)$ and $F(z)$ in $H$. 
Replacing $F$ by the composition
$p\circ F$, where $p:H\to H_0$ is the 
orthogonal projection, we may assume
that $H=H_0$ and the Hilbert space is at most 2-dimensional.
In particular the set of regular points is dense in $X$
by Corollary \ref{baer}. Hence we can assume
that $x,y$ and $z$ are regular points.
In particular
$C_y=C_y'\times H_y$ where $H_y$ is a Euclidean
space of dimension $\leq 2$ such that
$D_y F$ is the projection onto $H_y$.

Assume for a moment that
all $\tilde z$ near $z$ satisfy the equality 
$\tilde d (y,\tilde z) =0$. Then
$||F(\tilde z) -F(y)|| = d(\tilde z,y)$ and
since $F$ is 1-Lipschitz this implies
that all the initial vector of the geodesic
$y \tilde z$ lie in the $H_y$ factor of $C_y$.
It follows that $C_y'$ is trivial and
$F$ is an isometric embedding.
We are done in this case.

Hence   replacing $z$ by a nearby point
we may assume that  
$\tilde d(y,z) > 0$.  Set $\rho =\frac {d(y,z)}
{\tilde d (y,z)}$.

 Let $\gamma :[0,d(y,z)] \to X$ be a 
 unit speed geodesic between $y$
and $z$. Consider the function $h$ given by
$h(t)= \tilde d (x,y) + \tilde d(y, \gamma (t)) 
- \tilde d (x, \gamma (t))$.
We have $h(0)=0$, $h(d(y,z)) < 0$ and 
by equation (\ref{inpoint}) for $s>t$ the equality   
$$h(s)-h(t)= 
\tilde d (x, \gamma (t)) 
+ \tilde d(\gamma (t) ,\gamma (s)) - 
\tilde d (x,\gamma (s))$$ 
holds. 

The function $h$ is Lipschitz (hence differentiable 
almost everywhere)
and satisfies $h(t_0) <0$ for some $t_0$,
hence we find some $\varepsilon > 0$, $t \in [0,d(y,z)]$ such that
$h'(t)=-3\varepsilon$. 
Replacing $y$ by 
$\gamma (t)$ we may assume $t=0$. 
For a very small number $r<<\epsilon$ set
$z=\gamma (r)$. We then have $h(r) \leq -2\varepsilon r$.
Because of equation (\ref{inpoint}) we still have
$\rho =\frac {d(y,z)}{\tilde d (y,z)}$ and 
$\rho$ does not depend on
$\varepsilon$ and $r$.

Let $\eta _0:[0,s_0] \to X$ resp. $\eta _1: [0,s_1] \to X$ 
be geodesics
from $x$ to $y$ resp. to $z$. For 
$0<t\leq 1$ set $y_t = \eta _0 (t s_0)$
and $z_t =\eta _1 (t s_1)$.

We have $\tilde d(x,z_t)=t\tilde d(x,z)$; 
$\tilde d (x,y_t) = t\tilde
d(x,y)$. Moreover $||F(x)-F(z_t)|| = t||F(z)||$; 
$||F(x)-F(y_t)||=t ||F(y)||$
and $||F(z_t)-F(y_t)||=t||F(z)-F(y)||$. 

Since $X$ is a 
$CAT(\kappa)$ space and the pairwise 
distances of the points
$x,y,z$ are by assumption
$< \frac{\pi}{2 \sqrt{\kappa}}$,
there exists
$A \geq 0$ depending only on 
$\kappa$ such that
$$d(y_t,z_t) \leq t(d(y,z) + A\ d(y,z)^2)$$

We compute

\begin{align*}
\tilde d(y_t,z_t) &= 
\sqrt { d(y_t,z_t)^2 -t^2 ||F(y)-F(z)||^2} \\
&\leq \sqrt {t^2(d(y,z) +A\ d(y,z)^2)^2  
-t^2||F(y)-F(z)||^2}  \\  
&\leq t \sqrt {\tilde d (y,z) ^2  +B d(y,z) ^3} \\
&= t\tilde d(y,z) 
\sqrt {(1+  B \rho ^2 d(y,z))} \\
&\leq t \tilde d(y,z) (1 + C d(y,z))  
\end{align*}

for some constant $B$ depending only on
$\kappa$ and some constant $C$ depending 
only on $\rho$ and the curvature bound $\kappa$.
If $r = d(y,z)$ has been choosen small enough we 
thus obtain
\begin{equation} \label{est}
\tilde d(y_t,z_t) \leq t \tilde d(y,z) 
+ t \varepsilon \tilde d(y,z)
\end{equation}  
Since $h(r) \leq -2\varepsilon r$, $ r=d(y,z)$ we obtain
\begin{equation} \label{esta}
\tilde d (x,y) + \tilde d (y,z) - \tilde d (x,z)
\leq - 2 \varepsilon d(y,z) \leq - 2 \varepsilon  \tilde d (y,z)
\end{equation}
It follows that
\begin{align*}
\tilde d(x,z_t ) &=t \tilde d( x, z) \\
&\geq t(\tilde d (x,y) + \tilde d (y,z) 
+2\varepsilon \tilde d(y,z)) \\
& \geq \tilde d(x, y_t) + \tilde d(y_t,z_t) 
+ t \varepsilon \tilde d(y,z)
\end{align*}
where we used
equation (\ref{esta}) for the first
and equation (\ref{est}) for the second inequality.
Going to the limit $t\to 0$ we see that 
for the affine function 
$D_x F: C_x  \to H_x$ the correponding function 
$\tilde d_x : C_x \times C_x \to \R$,
$\tilde d_x(v,w) = \sqrt{d^2_x(v,w) - ||D_xF(v) - D_xF(w)||^2}$
is not a pseudo metric. 
But $D_x F: C_x \to  H_x$ is just the projection onto 
the Euclidean factor of $C_x = C_x' \times H_x$ 
(since $x$ is regular). Hence 
$\tilde d_x$  is just the metric on $C_x'$. 
Contradiction. 

This contradiction shows that
$\tilde d$ satisfies the triangle inequality
in the neighborhood of each point.
Using the $CAT(\kappa)$ property of $X$ and
equation (\ref{inpoint}) it is not difficult to
prove that the triangle inequality holds for
all triples of points.

.
\end{proof}

Let
$F:X \to H$ as in the assumption of Theorem 5.1,
then
$\tilde d$ defines a pseudometric on $X$.
Let $Y = X/\tilde d$ be the induced metric space.
A point in $Y$ is an equivalence class
$[x]$ where
$x\sim x'$ iff $\tilde d(x,x')=0$.
Theorem 5.1 implies immediately that the map
$X \to Y \times H$,
$x \mapsto ([x],F(x))$ is an isometric embedding.

For the proof of Theorem \ref{mainthm} we use
the affine map
$F= E_o: X \to \mathcal A^*$,
where $E_o$ is the evaluation map 
for some basepoint $o \in X$. 
By the discussion of section 3 and section 4,
the assumptions of Theorem 5.1 are satisfied.
Hence
$i:X \to Y \times \mathcal A^*$,
$x \mapsto ([x],E_o(x))$ is an isometric embedding.

If $f \in \tilde{\mathcal A}(X)$ is an affine
function on $X$, then define
$\hat{f}:\mathcal A^*  \to \R$ by
$$\hat{f}(\xi):= \xi ([f]) + f(o)$$
Then $\hat{f}$ is an affine function on $\mathcal A^*$
and
$$\hat{f}(E_o(x))=E_o(x)([f]) + f(o)=f(x)-f(o)+f(o)=f(x)$$
hence
$\hat{f}\circ p_{\mathcal A^*}\circ i = f$ as required.

We show now that $Y$ is a geodesic metric space.
Indeed if
$[y],[z] \in Y$
and
$\gamma:[0,1] \to X$
is a geodesic from
$x$ to $y$,
then
$t\mapsto [\gamma(t)]$ is a geodesic
in $Y$ due to equation (\ref{inpoint}).

We finally prove that the completion of
$Y$ is
$CAT(0)$ if $X$ is $CAT(0)$.
Let therefore
$[x],[y],[z] \in Y$ be arbitrary and
let
$[m] = [\gamma(\frac{1}{2})]$ be a midpoint
of $[y]$ and $[z]$.
We have to prove the Bruhat-Tits
$CAT(0)$ inequality (see e.g. \cite{BH} p.163):
\begin{equation} \label{estb}
\tilde d^2([x],[m]) \leq
\frac{1}{2} \tilde d^2([x],[y])
+ \frac{1}{2} \tilde d^2([x],[z])
- \frac{1}{4} \tilde d^2([y],[z])
\end{equation}
Since $X$ is $CAT(0)$ we have
$$ d^2(x,m) \leq
\frac{1}{2}  d^2(x,y)
+ \frac{1}{2}  d^2(x,z)
- \frac{1}{4}  d^2(y,z)$$
and since $F$ is affine we see
$$ ||F(x) - F(m)||^2 =
\frac{1}{2}  ||F(x) -F(y)||^2
+ \frac{1}{2}  ||F(x)-F(z)||^2
- \frac{1}{4}  ||F(y)-F(z)||^2.$$
Subtracting the two formulas we obtain inequality
(\ref{estb}).

\begin{tabbing}

Alexander Lytchak,\hskip10em\relax \= Viktor Schroeder,\\ 

Mathematisches Institut,\>
Institut f\"ur Mathematik, \\

Universit\"at Bonn,\> Universit\"at Z\"urich,\\
Beringstrasse 1, \>
 Winterthurer Strasse 190, \\

D-53115 Bonn, Germany\>  CH-8057 Z\"urich, Switzerland\\

{\tt lytschak@math.uni-bonn.de}\> 
{\tt vschroed@math.unizh.ch}\\
\end{tabbing}

\end{document}